\documentclass[]{article}
\usepackage{amssymb,amsmath,amscd,amsthm,xcolor}

\usepackage[a4paper]{geometry}
\usepackage{setspace}

\usepackage{authblk}
\usepackage{graphicx}
\usepackage{multirow}
\newcommand{\Keywords}[1]{\par\noindent
{\small{\em Keywords\/}: #1}}

\definecolor{pusu}{rgb}{0.0, 0.42, 0.24}

\newcommand{\Un}{\mathbf{1}}

\newcommand{\fy}{\varphi}

\newcommand{\tend}[1]{\longrightarrow}
\newcommand{\tends}[1]{\xrightarrow[#1]{}}
\newcommand{\tendsb}{\xrightarrow{a.s.}}
\newcommand{\tendsd}{\xrightarrow{\ d\ }}        
\newcommand{\norm}[1]{\lVert #1\rVert}                    

\newcommand{\abs}[1]{\lvert #1\rvert}                     

\DeclareMathOperator{\cov}{cov}              

\renewcommand{\d}{\mathrm{d}}

\let\0=\varnothing
\let\1=\Un

\newcommand{\edoc}{\end{document}}

\let\s=\sigma
\let\g=\gamma
\let\l=\lambda

\let\ls=\leqslant
\let\gs=\geqslant

\def\E{\hbox{\bf E}}

\let\kap=\varkappa

\def\dvit{\colon\ }
\def\DD{\mathop{\bf Var}\nolimits}
\def\cov{\mathop{\rm cov}\nolimits}
\def\diag{\hbox{\rm diag}}

\theoremstyle{plain}
\newtheorem{thm}{Theorem}

\newtheorem{prop}[thm]{Proposition}
\newtheorem{dfn}[thm]{Definition}
\newtheorem{rem}[thm]{Remark}
\newtheorem{cor}[thm]{Corollary}

\numberwithin{equation}{section}

\title{CLT for quadratic variation of Gaussian processes and its application to the estimation of the Orey index.}
\date{}
\author{K. Kubilius}

\begin{document}

\maketitle

\vskip-0.5cm

\centerline{\small Vilnius University, Faculty of Mathematics
and Informatics} \centerline{\small Akademijos 4, LT-08663 Vilnius,
Lithuania}


\abstract{We give a  two-dimensional central limit theorem (CLT) for the second-order quadratic variation of the centered Gaussian processes on $[0,T]$. Though the approach we use is well known in the literature, the conditions under which the CLT holds are usually based on differentiability of the corresponding covariance function. In our case, we replace differentiability conditions by the convergence of the scaled sums of the second-order moments. To illustrate the usefulness and easiness of use  of the approach, we apply the obtained CLT  to proving  the asymptotic normality of the estimator of the Orey index of a subfractional Brownian motion.

\bigskip
 }

\vskip1cm

\section{Introduction}\label{s:intro}
In recent decades, self-similar Gaussian processes attracted vast attention due to the burst of many successful applications in different areas including (but not limited to) financial sector, telecommunications, hydrology, biology, etc. Consequently, many theoretical investigations of the properties of these processes were made. Since the behavior of these processes crucially depends on the value of the corresponding self-similarity index $H$, the very important task from a statistical point of view is an estimation  of the value of $H$ having sample data. As a rule, the latter task involves proving asymptotic normality of the constructed estimator, and, in different scenarios, different CLTs may appear useful. This was one of the major reasons for presenting CLT given in the sequel as this theorem targets the case of a discrete single trajectory sample data $\{X_{\frac{1}{n}T},\dots,X_{\frac{n-1}{n}T},X_{\frac{n}{n}T}\}$ of the observed Gaussian process $(X_t)_{t\in[0,T]}$. Another reason stems from the fact that there were various attempts to embed self-similar Gaussian processes into the larger classes and then treat them as the separate cases of these classes. In particular, Orey \cite{or} was among the first who proposed one of such classes, later considered in a more general form by Norvai\v sa \cite{rn0}, \cite{rn1}, \cite{rn2} and Kubilius \cite{kk1} (see also \cite{KMR}), who gave generalized definition of the Orey class. In \cite{rn2} and \cite{kk1}, it was demonstrated that this class includes such popular models as fractional Brownian motion (fBm, \cite{Kolmogorov40}, \cite{Mandelbrot68}), bi-fractional Brownian motion (bifBm, \cite{HV}) and sub-fractional Brownian motion (sfBm, \cite{BGT}) among the rest. Moreover, in \cite{kk1}, there were proposed consistent estimators of the Orey index under the sampling setting described above. However, asymptotic normality of the suggested estimators was left open. The present paper aims  to fill this gap.

The given CLT is based on the second-order quadratic variations.
The literature devoted to the behaviour of second-order quadratic variations and CLTs theorems of this type in particular is very abundant (e.g., Bardet and Surgailis \cite{BS11, BS13}, Baxter \cite{baxter}, Gladyshev \cite{glad}, Klein and Gin\'e \cite{kg}, Guyon and Le\'on \cite{guyon}, Istas and Lang \cite{IL}, Benassi et. al. \cite{Istas98}, Coeurjolly \cite{Coeurjolly-01}, Cohen et.al. \cite{guyon}, B\'egyn \cite{begyn1,begyn2}, Norvai\v sa \cite{rn1,rn2}, Malukas \cite{ma}, Kubilius \cite{kk1}, Vitasaari \cite{viit}).

For the asymptotic normality of the proposed estimators in \cite{kk1} it is very natural to use the CLT for a second-order quadratic variation of Gaussian processes  obtained by Begyn in \cite{begyn2}. Unfortunately, his theorem is not always applicable to the class of processes considered by us. To see this, consider, for instance, sfBm and bifBm (see Appendix).

The paper is organized as follows. In Section \ref{s:notation}, we introduce notions and, following \cite{kk1}, restate the definition of the Gaussian process belonging to the class of Gaussian processes having Orey index. In Section \ref{s:mainRes}, we state the main theorem and several direct corollaries. In Section \ref{s:applications} we demonstrate that the theorem can be very effective for derivation of the asymptotic normality of estimators of the Orey index for sfBm. Section \ref{s:proofs} contains proofs. At the end of the paper there is an Appendix on the application of Theorem 2 of  \cite{begyn2} for sfBm and bifBm.

\section{Notation and auxiliary results}\label{s:notation}

Let $X = (X_t)_{t \in [0, T ]}$ be a second-order stochastic process with an
incremental variance function
\[
\s_X^2(s, t) := \E[X_t-X_s]^2,\quad (s, t)\in[0, T ]^2,
\]
and let $RV_\rho$ be the class of functions defined on $(0,T)$ and regularly varying at zero with an index of variation equal to $\rho\in\mathbb{R}$. Define
\[
    \Psi =\left\{f\in RV_{1}\ \Big\vert\ L(h)=\frac{f(h)}{h}\tends{h\to 0+0}\infty\right\},
\]
and for any $\fy\in\Psi$,

\begin{align*}
\g_*(\varphi):=&\inf\bigg\{\g>0\dvit \lim_{h\downarrow 0}\sup_{\varphi(h)\ls s\ls T-h}\frac{h^\g}{\s_X(s,s+h)}=0\bigg\},\\
\widetilde\g_*:=&\inf\Big\{\g>0\dvit \lim_{h\downarrow 0}\frac{h^\g}{\s_X(0,h)}=0\Big\}
\end{align*}
and
\begin{align*}
\g^*(\varphi):=&\sup\bigg\{\g>0\dvit \lim_{h\downarrow 0}\inf_{\varphi(h)\ls s\ls T-h}\frac{h^\g}{\s_X(s,s+h)}=+\infty\bigg\},\\
\widetilde\g^*:=&\sup\Big\{\g>0\dvit \lim_{h\downarrow 0}\frac{h^\g}{\s_X(0,h)}=+\infty\Big\}\,,
\end{align*}
where $\varphi\in\Psi$. Note that $0\ls\widetilde\g^*\ls\widetilde\g_*\ls +\infty$ and $0\ls\g^*(\varphi)\ls\g_*(\varphi)\ls +\infty$.

\begin{dfn}[\cite{kk1}, see also \cite{KMR}]\label{oreydef1}
Assume that $\sup_{0\ls s\ls T-h}\s_X(s,s+h)\to 0$ as $h\to 0$. If $\g_*(\varphi)=\widetilde\g_*=\g^*(\varphi)=\widetilde\g^*$ for any function $\varphi\in\Psi$, then we say that the process $X=(X_t)_{t\in[0,T]}$ has the Orey index\index{Orey index!extended} $\g_X=\widetilde\g_*=\widetilde\g^*$.
\end{dfn}

\begin{rem}
In case when $X$ has stationary increments, one only needs to check that $\widetilde\g_*=\widetilde\g^*$.
\end{rem}
 In what follows, we will make use of the following result.

\begin{thm}[\cite{kk1}, see also \cite{KMR}]\label{th:Orey index} Assume that, for some $\g\in(0,1)$, the zero-mean second-order stochastic process $X=(X_t)_{t\in[0,T]}$ satisfies conditions:

(C1)\quad $\s_X(0,\delta)\asymp\delta^{\g}$, i.e., $\s_X(0,\delta)$ and $\delta^{\g}$ is of the same order as $\delta\downarrow 0$;

(C2)\quad there exists a constant $\kap>0$ such that
\[
 \L(\delta):=\sup_{\varphi(\delta)\ls t\ls T-\delta} \sup_{0<h\ls \delta}\bigg\vert\frac{\s_X(t,t+h)}{\kap h^{\g}}- 1\bigg\vert\longrightarrow 0\qquad\mbox{as}\ \delta\downarrow 0
\]
for every function $\varphi\in\Psi$.

Then it has the Orey index equal to $\g$.
\end{thm}
Let $X=(X_t)_{t\in[0,T]}$, $T>0$, $X_0=0$, be a zero-mean Gaussian process. Finishing the section, we introduce some notions repeatedly used in the sequel:

\begin{align}\label{e:deltas}
V^{\widehat X}_{in,T} =&\sum_{k=1}^{in-1} \big(\Delta^{(2)}_{in,k}\widehat X\big)^2, \qquad i=1,2, \\
d_{j,k}^{\widehat X,in}=&\E \Delta^{(2)}_{in,j}\widehat X\Delta^{(2)}_{in,k}\widehat X,\qquad 1\ls j,k\ls in-1,\quad i=1,2,\\
c_{j,k}^{\widehat X,n}=&\E\Delta^{(2)}_{n,j}\widehat X \Delta^{(2)}_{2n,k}\widehat X,\quad 1\ls j\ls n-1,\  1\ls k\ls 2n-1,\\
\Delta^{(2)}_{in,k} X =&X_{\frac{k+1}{in}T}-2 X_{\frac{k}{in}T}+ X_{\frac{k-1}{in}T}\quad 1\ls k\ls in-1,\quad i=1,2,
\end{align}
where $\Delta^{(2)}_{in,k}\widehat X$ are normalized second-order increments defined by
\[
\Delta^{(2)}_{in,k}\widehat X=
\begin{cases}
\displaystyle{\frac{\Delta^{(2)}_{in,k}X}{\sqrt{\E (\Delta^{(2)}_{in,k}X)^2}}} & \mbox{for Gaussian process},\\
\displaystyle{\frac{\Delta^{(2)}_{in,k}X}{\kap\sqrt{4-2^{2\g}}\,(\frac{T}{in})^\g}} & \mbox{for Gaussian process with the Orey index}\ \g,
\end{cases}
\]
with $\kap$ equal to normalizing constant in condition (C2) above.
Finally, by $B^\g=(B^\g_t)_{t\in[0,T]}$ we denote the fBm having Hurst index $\gamma$, which is also its Orey index. Recall that the covariance function of $B^\gamma$ is given by
\begin{equation*}
    \E B^{\gamma}_tB^\gamma_s=\frac{1}{2} \left(t^{2\gamma}+s^{2\gamma}-\abs{t-s}^{2\gamma}\right),
\end{equation*}
and that $B^\gamma$ is the only (up to the constant multiplier) self-similar Gaussian process having stationary increments.

\section{Main results}\label{s:mainRes}

The theorem below is the main result of the paper.

\begin{thm}\label{t:main} Let $X=(X_t)_{t\in[0,T]}$, $T>0$, be a zero-mean Gaussian process, $X_0=0$.  Assume that, for $\forall n$,
\begin{equation}\label{e:maks_sum}
\max_{1\ls k\ls in-1}\sum_{j=1}^{in-1}\vert d_{j,k}^{\widehat X,in}\vert \ls C,\qquad i=1,2,
\end{equation}
for some constant $C$  and
\begin{equation}\label{e:kovariaciju_riba}
n \cov\big((in)^{-1}V^{\widehat X}_{in,T},(jn)^{-1}V^{\widehat X}_{jn,T}\big) \longrightarrow \Sigma_{ij}\qquad\mbox{for all}\quad i,j\in\{1,2\}.
\end{equation}
Then
\begin{equation}\label{e:lim1}
\sqrt{n}\left(
      \begin{array}{c}
       n^{-1} \big(V^{\widehat X}_{n,T}-\E V^{\widehat X}_{n,T}\big) \\
        (2n)^{-1}\big(V^{\widehat X}_{2n,T}-\E V^{\widehat X}_{2n,T}\big) \\
      \end{array}
    \right)\tendsd \mathcal{N}\left(0;\Sigma\right),\qquad
    \Sigma= \left(\begin{array}{cc}
                 \Sigma_{11} & \Sigma_{12} \\
                 \Sigma_{12} & \Sigma_{22}  \\
               \end{array}
             \right),
\end{equation}
where $\mathcal{N}\left(0;\Sigma\right)$ is a Gaussian vector with the covariance matrix $\Sigma$.
\end{thm}

Note that, by \eqref{e:deltas}, $\E V^{\widehat X}_{in,T}=in-1$, therefore \eqref{e:lim1} can be replaced by
\begin{equation}\label{e:lim1modif}
    \sqrt{n}\left(
      \begin{array}{c}
       n^{-1} \big(V^{\widehat X}_{n,T}-1\big) \\
        (2n)^{-1}\big(V^{\widehat X}_{2n,T}-1\big) \\
      \end{array}
    \right)\tendsd \mathcal{N}\left(0;\Sigma\right).
\end{equation}
Another important observation is related to the case of the fBm $B^{\gamma}$. For that particular case, \eqref{e:lim1} was shown to hold by several authors (see, e.g., \cite{cohen}, \cite{begyn2}). Hence, for this case, $\Sigma=\Sigma_\gamma$, where
\begin{align}
\Sigma_{11}=&2\bigg(1+\frac{2}{(4-2^{2\g})^2}\sum_{j=1}^{\infty}\widehat\rho_\g^{\,2}(j)\bigg), \qquad \Sigma_{22}=\frac 12\, \Sigma_{11},\nonumber\\
\Sigma_{12}=&\Sigma_{21}=\frac{1}{2^{2\g}(4-2^{2\g})^2}\sum_{j\in\mathbb{Z}} \widetilde\rho_\g^{\,2}(j),\nonumber\\
\widehat\rho_\g(j)=&\frac 12\,\big[-6\abs{j}^{2\g} -\abs{j-2}^{2\g} -\abs{j+2}^{2\g}+4\abs{j-1}^{2\g} +4\abs{j+1}^{2\g}\big],\label{e:covariance1}\\
\widetilde{\rho}_\g(j)=&\frac{1}{2}\big[\abs{j+1}^{2\g} +2\abs{j+2}^{2\g}-\abs{j+3}^{2\g} + \abs{j-1}^{2\g}-4\abs{j}^{2\g}\nonumber\\
&- \abs{j-3}^{2\g} + 2\abs{j-2}^{2\g}\Big].\nonumber
\end{align}
Note that for $\widehat\rho_\g(j)$ and $\widetilde{\rho}_\g(j)$ we have estimates (see \cite{KMR} 63 p., 65 p.)
\begin{align}
\frac{\vert \widehat\rho_\g(j)\vert}{4-2^{2\g}}\ls& \frac{243}{20\ln 4}\,j^{2\g-4}\ls 9 k^{2\g-4}\qquad \mbox{for any}\ j\gs 3, \label{e:vanishing0}\\
\frac{\vert \widetilde \rho_H(j)\vert}{4-2^{2\g}}\ls& \frac{36}{\ln 4}\, j^{2H-4}\ls 26 j^{2H-4}\qquad \mbox{for any}\ j\gs 4.\label{e:vanishing}
\end{align}
We make use of that fact in Section \ref{s:proofs}.

Combining Theorem \ref{t:main} and \eqref{e:covariance1} leads to the following corollaries.

\begin{cor}\label{c:cor1} Let $X=(X_t)_{t\in[0,T]}$, $T>0$, $X_0=0$, be a zero-mean Gaussian process. Assume that condition \eqref{e:maks_sum} holds and there exists a fBm  $B^\g$  such that
\begin{equation}\label{e:covLim}
    \frac{1}{n}\left(\cov\left(i^{-1}V_{in,T}^{\widehat X},j^{-1}V_{jn,T}^{\widehat X}\right)- \cov\left(i^{-1}V_{in,T}^{\widehat B^{\gamma}},j^{-1}V_{jn,T}^{\widehat B^{\gamma}}\right)\right)\tends{n\to\infty}0
\end{equation}
for all $i,j\in\{1,2\}$. Then \eqref{e:lim1} applies with $\Sigma$ equal to $\Sigma_\gamma$ given by \eqref{e:covariance1}.
\end{cor}

\begin{rem}\label{r:sufficientCond}
Note that for $i=1,2$
\begin{align}
&\frac 1n\sum_{k=1}^{in-1}\big\vert (d^{\widehat X,in}_{k,k})^2 -(d^{\widehat B^\g,in}_{k,k})^2\big\vert+\frac 1n\sum_{k=1}^{in-2}\sum_{m=1}^{in-k-1}\big\vert (d^{\widehat X,in}_{k,k+m})^2 -(d^{\widehat B^\g,in}_{k,k+m})^2\big\vert  \longrightarrow 0,\label{e:kovariaciju_skirtumas1}\\
&\frac 1n\sum_{j=1}^{n-1}\sum_{k=1}^{2n-1} \big\vert (c_{j,k}^{\widehat X,n})^2 - (c_{j,k}^{\widehat B^\g,n})^2\big\vert\longrightarrow 0\label{e:kovariaciju_skirtumas2}
\end{align}
is a sufficient condition for \eqref{e:covLim} to hold. This follows from the equalities
\begin{align}
\DD \sum_{k=1}^{in-1} \left(\Delta^{(2)}_{in,k}\widehat X\right)=&2\sum_{k=1}^{in-1} (d^{\widehat X,in}_{kk})^2 +4\sum_{k=1}^{in-2}\sum_{m=1}^{in-k-1}(d^{\widehat X,in}_{k,k+m})^2,\label{e:issert1}\\
\cov(V^{\widehat X}_{n,T},V^{\widehat X}_{2n,T}) =& 2\sum_{j=1}^{n-1}\sum_{k=1}^{2n-1}\big(\E\Delta^{(2)}_{n,j}\widehat X \Delta^{(2)}_{2n,k}\widehat X\big)^2\label{e:issert2}
\end{align}
which are proved using the Isserlis formulas (see \cite{Iss}).
\end{rem}

\begin{cor}\label{c:cor2} Assume that $X=(X_t)_{t\in[0,T]}$, $T>0$, $X_0=0$, is a zero-mean Gaussian process having the Orey index $\g$.  Also, assume that conditions stated in Corollary \ref{c:cor1} hold and
\begin{equation}\label{e:dispersiju_skirtumas3}
\sqrt{n}\big((in)^{-1}\E V^{\widehat X}_{in,T}-1\big)\tend{} 0.
\end{equation}
Then \eqref{e:lim1modif} applies with $\Sigma$ equal to $\Sigma_\gamma$ given by \eqref{e:covariance1}.
\end{cor}

\begin{rem}
Note that, by the definition of $V^{\widehat X}_{in,T},\Delta^{(2)}_{in,k}\widehat X$, in Corollary \ref{c:cor2}, we normalize $\Delta^{(2)}_{in,k} X$ by $\kap\sqrt{4-2^{2\gamma}}\left(\frac{T}{in}\right)^{\gamma}$ rather than $\E(\Delta^{(2)}_{in,k} X)^2$. That is, in this case, we do not have an equality $\E V^{\widehat X}_{in,T}=in-1$ and \eqref{e:dispersiju_skirtumas3} needs verification.
\end{rem}

We finish this section by defining an estimator of the Orey index and providing a theorem which illustrates an application of the main Theorem \ref{t:main} to the estimation problem announced in the introductory Section \ref{s:intro}. The given estimator was previously considered for the case of the fBm and some other processes as well by \cite{Istas98}, \cite{begyn2}. It also appeared to be consistent for the case of a more general class of Gaussian processes having the Orey index (see \cite{kk1} and \cite{KMR}) and for solutions of the stochastic differential equations driven by the fBm (see \cite{kubmel1}, \cite{ksm}, and \cite{KMR}).

\begin{thm}\label{t:gNormality} Let $X=(X_t)_{t\in[0,T]}$, $T>0$, $X_0=0$, be a Gaussian process satisfying assumptions of Corollary \ref{c:cor2}. Define
\[
\widehat{\g}_n=\frac{1}{2}-\frac{1}{2\ln
    2}\ln\left(\frac{\sum_{k=1}^{2n-1} \left(\Delta^{(2)}_{2n,k}X\right)^2}{\sum_{k=1}^{n-1} \left(\Delta^{(2)}_{n,k}X\right)^2}\right).
\]
Then
\[
\widehat{\g}_n\to \gamma\quad\mbox{a.s.}
\]
and
\[
2\ln 2\,\sqrt{n}\,(\widehat{\g}_n-\g)\tendsd \mathcal{N}(0;\sigma_\g^2),\qquad \sigma_\g^2=\frac 32\, \Sigma_{11}-2 \Sigma_{12}.
\]
\end{thm}

\section{Exemplary application}\label{s:applications}

In this section, we demonstrate an application of the main result for the case of sfBm. Recall that it is a centered Gaussian process having covariance function
\begin{equation}
    G_H(s,t):=s^{2H} +t^{2H}-\frac{1}{2}\big[(s+t)^{2H}+\vert s-t\vert^{2H}\big],\quad (s,t)\in[0,T]^2.
\end{equation}
For this process, the Orey index $\gamma$ is equal to the value of the process defining parameter $H\in(0,1)$, and the normalizing constant $\kap$ in (C2) of Theorem \ref{th:Orey index} is equal to 1.

\begin{prop}\label{p:examples}
Let $S^H=(S^H_t)_{t\in[0,T]}$ be a sfBms. Then the Theorem \ref{t:gNormality} applies to $S^H$.
\end{prop}
As noted in the introduction, we cannot apply B\'egyn's result from \cite{begyn2} for sfBm here. Using our proven theorem we get a simple proof that leads very quickly to the desired result.

\section{Proofs}\label{s:proofs}
In this section, when writing bounds, it always suffices just to know that the left hand side is bounded by some constant independent of $n$ tending to $\infty$. Therefore, to suppress notation, $C$ always stands for such a constant. It may vary from line to line and depend on $\gamma$ and/or other fixed quantities but not on $n\to\infty$. The fact that $C$ has changed is in no way indicated.

\smallskip\noindent{\bf Proof of Theorem \ref{t:main}.} Let
\[
\mathbf{X}_n=\sqrt{n}\left(
      \begin{array}{c}
       n^{-1} \big(V^{\widehat X}_{n,T}-\E V^{\widehat X}_{n,T}\big) \\
        (2n)^{-1}\big(V^{\widehat X}_{2n,T}-\E V^{\widehat X}_{2n,T}\big) \\
      \end{array}
    \right).
\]
To find out limiting distribution of $\mathbf{X}_n$, we compute
a limiting moment generating function $\lim_{n\to\infty} M_{\mathbf{X}_n}(\mu)=M(\mu)$, where $\mu=(\mu_1,\mu_2)$.

Define a centered Gaussian vector $\mathbf{G}_n=(G_n^{(j)}, 1\ls j\ls 3n-2)$ as follows:
\begin{align*}
G_n^{(j)}=& \Delta^{(2)}_{n,j} \widehat X,\quad 1\ls j\ls n-1,\\
G_n^{(j)}=& \sqrt{2^{-1}}\,\Delta^{(2)}_{2n,j+1-n} \widehat X,\quad n\ls j\ls 3n-2
\end{align*}
Let $\Sigma_{\mathbf{G}_n}$ be the covariance matrix of $\mathbf{G}_n$. Set
$\widetilde D_n =\left (\Sigma^{1/2}_{\mathbf{G}_n}\right )^T D_n \Sigma^{1/2}_{\mathbf{G}_n}$, where $D_n$ is a diagonal matrix given by
\[
D_n=\diag(\underbrace{\mu_1,\ldots,\mu_1}_{n-1}, \underbrace{\mu_2,\ldots,\mu_2}_{2n-1}).
\]
For a square matrix $A$, denote its eigenvalues by $\lambda_k$, maximal value among $\abs{\lambda_k}$'s by $\rho(A)$, and the operator norm $\sup_{\norm{x}=1}\norm{A x}$ by $\norm{A}$. For symmetric matrix $\widetilde D_n$ its norm is equal to its spectral norm, i.e. $\norm{\widetilde D_n}=\rho(\widetilde D_n) :=\max_{k}\vert\lambda_{k}(\widetilde D_n)\vert$. Since norm $\norm{\,\cdot}$ is submultiplicative norm then
\begin{align*}
\rho(\widetilde D_n)=&\norm{\widetilde D_n} \ls \norm{\Sigma^{1/2}_{\mathbf{G}_n}}\cdot \norm{D_n}\cdot \norm{\Sigma^{1/2}_{\mathbf{G}_n}}=\norm{\Sigma^{1/2}_{\mathbf{G}_n}}^2\cdot \norm{D_n}
=\rho^2\big(\Sigma^{1/2}_{\mathbf{G}_n}\big)\cdot\rho(D_n)\\ =&\rho(\Sigma_{\mathbf{G}_n})\cdot\max\{\abs{\mu_1}, \abs{\mu_2}\}.
\end{align*}
Recall that the largest eigenvalue of a symmetric non-negative definite matrix does not exceed maximal row (col) sum of absolute values. Thus,
\begin{equation}\label{i:norm}
    \rho(\Sigma_{\mathbf{G}_n}) \ls\max_{j}\sum_{i=1}^{3n-2}\big\vert(\Sigma_{\mathbf{G}_n})_{ij}\big\vert\,.
\end{equation}
Next, note that
\begin{align*}
(\Sigma_{\mathbf{G}_n})_{i,j+1-n}=&\frac{1}{\sqrt 2}\,\E \big[ \Delta^{(2)}_{n,i}\widehat X\Delta^{(2)}_{2n,j+1-n} \widehat X\big]\\
=&\frac{1}{\sqrt 2}\,\E\big[ \big( \Delta^{(2)}_{2n,2i+1}\widehat X+\Delta^{(2)}_{2n,2i-1}\widehat X +2\Delta^{(2)}_{2n,2i}\widehat X\big)\Delta^{(2)}_{2n,j+1-n}\widehat X\big]\nonumber\\
=&\frac{1}{\sqrt 2}\,\big[d_{2i+1,j+1-n}^{\widehat X,2n}+ d_{2i-1,j+1-n}^{\widehat X,2n}+2d_{2i,j+1-n}^{\widehat X,2n}\big]
\end{align*}
for $1\ls i\ls n-1$ and $n\ls j\ls 3n-2$. Therefore, by  (\ref{i:norm}) and  (\ref{e:maks_sum}),
\begin{equation*}
\rho(\Sigma_{\mathbf{G}_n})\ls C\left( \max_{1\ls j\ls n-1}\sum_{i=1}^{n-1}\vert d_{i,j}^{\widehat X,n}\vert +\,\max_{1\ls j\ls 2n-1}\sum_{i=1}^{2n-1}\vert d_{i,j}^{\widehat X,2n}\vert \right)\ls C.
\end{equation*}
Summing up, we come to conclusion that $\rho(\widetilde D_n)$ is uniformly (in $n$, $k$) bounded by finite constant depending only on $\mu_1$, $\mu_2$.


Recall that $\mathbf{G}_n\stackrel{d}{=}\sqrt{\Sigma_{\mathbf{G}_n}}\,\mathbf{Z}_n$ with $\mathbf{Z}_n\sim \mathcal{N}(0;I_{3n-2})$, where  $\stackrel{d}{=}$ denotes equality in distribution, $I_{3n-2}$ denotes an identity $(3n-2)$ matrix. Therefore,
\begin{align}\label{e:orto}
\mathbf{G}_n^T D_n \mathbf{G}_n\stackrel{d}{=}&(\Sigma^{1/2}_{\mathbf{G}_n} \mathbf{Z}_n)^T D_n \Sigma^{1/2}_{\mathbf{G}_n} \mathbf{Z}_n =\mathbf{Z}_n^T(\Sigma^{1/2}_{\mathbf{G}_n})^T D_n \Sigma^{1/2}_{\mathbf{G}_n} \mathbf{Z}_n=\mathbf{Z}_n^T\widetilde D_n \mathbf{Z}_n.
\end{align}
Let $\widetilde D_n=Q_n^T\Lambda(\widetilde D_n)Q_n$ be canonical representation of $\widetilde D_n$
via diagonal matrix of eigenvalues and corresponding orthogonal
matrix of eigenvectors. Since orthogonal transform does not change the distribution of $\mathbf{Z}_n$,
\begin{equation}\label{e:trace}
\mathbf{Z}_n^T\widetilde D_n \mathbf{Z}_n= \mathbf{Z}_n^TQ_n^T\Lambda(\widetilde D_n)Q_n \mathbf{Z}_n \stackrel{d}{=}\mathbf{Z}_n^T\Lambda(\widetilde D_n) \mathbf{Z}_n=\sum_{j=1}^{3n-2}Z_{n,j}^2\l_j(\widetilde D_n).
\end{equation}
By the above, $\rho(\widetilde D_n)\ls C$. Therefore, for $n$ sufficiently large, $n^{-1/2}\max_k\vert\l_k(\widetilde D_n)\vert <1/2$, and, in what follows, we assume without loss of generality that this condition holds for all $n$. Under this assumption, $M_{\mathbf{X}_n}(\mu)$ is well defined for all $n$ and since
\[
\mu^T\mathbf{X}_n=\sqrt{n}(\mu_1,\mu_2)\left(\begin{array}{c}
       n^{-1} \big(V^{\widehat X}_{n,T}-\E V^{\widehat X}_{n,T}\big)\\
       (2n)^{-1} \big(V^{\widehat X}_{2n,T}-\E V^{\widehat X}_{2n,T}\big) \\
      \end{array}
    \right)=\frac{1}{\sqrt{n}}\big(\mathbf{G}_n^T D_n \mathbf{G}_n-\E \mathbf{G}_n^T D_n \mathbf{G}_n\big),
\]
\eqref{e:orto}--\eqref{e:trace} implies
\begin{align*}
M_{\mathbf{X}_n}(\mu)=&\exp\Big\{-\frac{\E\, \mathbf{G}_n^T D_n \mathbf{G}_n}{\sqrt{n}}\Big\}\E\exp\bigg\{\sum_{j=1}^{3n-2}Z_{n,j}^2 \frac{\l_j(\widetilde D_n)}{\sqrt n}\bigg\},
\end{align*}
where $Z_{n,j}^2$ are i.i.d. and each $Z_{n,j}^2$ has chi-square distribution $\chi^2(1)$ with $M_{Z_{n,j}^2}(x)=\frac{1}{\sqrt{1-2x}}$ for $ x\in(-1/2,1/2)$. Thus,
\begin{align*}
M_{\mathbf{X}_n}(\mu)=&\exp{\left\{-\frac{\E\, \mathbf{G}_n^T D_n \mathbf{G}_n}{\sqrt{n}}\right\}}\prod_{j=1}^{3n-2}M_{\chi^2(1)}
\left(\frac{\l_j(\widetilde D_n)}{\sqrt{n}}\right)\\
=&\exp{\left\{-\frac{\E\, \mathbf{G}_n^T D_n \mathbf{G}_n}{\sqrt{n}}\right\}}\left(\prod_{j=1}^{3n-2}\frac{1}{1 -2\,\frac{\l_j(\widetilde D_n)}{\sqrt{n}}}\right)^\frac{1}{2}\\
=&\exp{\left\{-\frac{\E\, \mathbf{G}_n^T D_n \mathbf{G}_n}{\sqrt{n}}-\frac{1}{2}\sum_{j=1}^{3n-2} \ln\left(1-2\,\frac{\l_j(\widetilde D_n)}{\sqrt{n}}\right)\right\}}\\
=&\exp{\left\{-\frac{\E \mathbf{G}_n^T D_n \mathbf{G}_n}{\sqrt{n}}+\frac{1}{2} \sum_{j=1}^{3n-2}\left(2\,\frac{\l_j(\widetilde D_n)}{\sqrt{n}} +4\,\frac{\l_j^2(\widetilde D_n)}{2{n}}\right) +O\left(\frac{1}{\sqrt{n}}\right)\right\}}\\
    =&\exp{\left\{\frac{1}{n}\sum_{j=1}^{3n-2} \l_j^2(\widetilde D_n)\right\}}\exp{\left\{O\left(\frac{1}{\sqrt{n}}\right)\right\}},
\end{align*}
where we have expanded $x\mapsto \ln(1-x)$ in a neighborhood of 0 and used \eqref{e:orto}--\eqref{e:trace} to deduce that $\E\mathbf{G}_n^T D_n \mathbf{G}_n=\sum_{j=1}^{3n-2}\l_j(\widetilde D_n)$. Therefore, it remains to compute limiting value of the first multiplier. By the definition of $\widetilde{D}_n$,
\begin{align}\label{e:expr_for_trace}
    \sum_{j=1}^{3n-2}\l_j^2(\widetilde D_n)=&\mathrm{tr}(\widetilde{D}_n^2)=
    \mathrm{tr}\big(((\sqrt{\Sigma_{\mathbf{G}_n}}\,)^T D_n\sqrt{\Sigma_{\mathbf{G}_n}}\,)^2\big)=
    \mathrm{tr}\big((D_n\Sigma_{\mathbf{G}_n})^2\big)\nonumber\\
    =&\sum_{i=1}^{3n-2}\sum_{j=1}^{3n-2}(D_n\Sigma_{\mathbf{G}_n})_{ij}(D_n \Sigma_{\mathbf{G}_n})_{ji}.
\end{align}
Rearranging \eqref{e:expr_for_trace}  yields
\begin{align}\label{e:squareTracexpanded}
&\sum_{i=1}^{3n-2}\sum_{j=1}^{3n-2}(D_n\Sigma_{\mathbf{G}_n})_{ij}(D_n \Sigma_{\mathbf{G}_n})_{ji} \nonumber\\
&\quad=\mu_1^2\sum_{i=1}^{n-1}\sum_{j=1}^{n-1} d_{i,j}^{\widehat X,n} d_{j,i}^{\widehat X,n} +\mu_1 \mu_2\sum_{i=1}^{n-1}\sum_{j=1}^{2n-1} c_{i,j}^{\widehat X,n}c_{j,i}^{\widehat X,n} +\frac{\mu_2^2}{4} \sum_{i=1}^{2n-1}\sum_{j=1}^{2n-1} d_{i,j}^{\widehat X,n}d_{j,i}^{\widehat X,2n}.
\end{align}
By applying equalities \eqref{e:issert1} and \eqref{e:issert2} we get
\begin{align*}
&\sum_{i=1}^{3n-2}\sum_{j=1}^{3n-2}(D_n\Sigma_{\mathbf{G}_n})_{ij}(D_n \Sigma_{\mathbf{G}_n})_{ji}\\
&\quad=\frac{1}{2}\left(\mu_1^2\DD V_{n,T}^{\widehat X} +2\mu_1\mu_2\cov\left(V^{\widehat X}_{n,T},V^{\widehat X}_{2n,T}\right) +\frac{\mu_2^2}{4}\DD \left(V_{2n,T}^{\widehat X}\right)\right).
\end{align*}
Taking into account assumption \eqref{e:kovariaciju_riba}, dividing each sum by $n$, and passing to the limit yields $M(\mu)=\exp\{\frac 12\mu^T\Sigma \mu\}$.
\endproof

\smallskip{\bf Proof of Corollaries \ref{c:cor1}, \ref{c:cor2}.} The statements are obvious.

\smallskip{\bf Proof of Theorem \ref{t:gNormality}.} Under conditions of the theorem we have (see \cite{kk1}, Corollary 3.10, and \cite{KMR}, Corollary 6.13),
\begin{equation}\label{e:convergence}
\bigg(\frac{n}{T}\bigg)^{2\g-1}\sum_{k=1}^{n-1} \big(\Delta^{(2)}_{n,k}X\big)^2\tendsb{}\kap^2(4-2^{2\g})T.
\end{equation}
Strong consistency and asymptotic normality of the estimator $\widehat{\gamma}_n$ follows from the expression
\begin{align*}
\widehat{\g}_n=&\frac 12-\frac{1}{2\ln
    2}\ln\left(\frac{2^{1-2\g}\frac{1}{2n}V^{\widehat X}_{2n,T}}{\frac{1}{n}V^{\widehat X}_{n,T}}\right)
=\g-\frac{1}{2\ln 2}\ln\left(\frac{\frac{1}{2n}V^{\widehat X}_{2n,T}}{\frac{1}{n}V^{\widehat X}_{n,T}}\right).
\end{align*}
Really, from \eqref{e:convergence} the strong consistency holds. Corollary \ref{c:cor2} and the Delta method yield (see Remark 2.12 \cite{KMR})
\[
2\ln 2\sqrt{n}\,\big(\widehat{\g}_n-\g\big)= \sqrt n \left(\ln\bigg(\frac{\frac{1}{2n}V^{\widehat X}_{2n,T}}{\frac{1}{n}V^{\widehat X}_{n,T}}\bigg)-\ln 1\right)\tendsd \mathcal{N}\left(0;\s^2_\g\right).
\]
\endproof

\smallskip{\bf Proof of Proposition \ref{p:examples}.} Direct calculations show that, for $0\ls u<v\ls s<t$,
\begin{align}\label{e:difference}
&\E\big[\big(S^H_v-S^H_u\big)\big(S^H_t-S^H_s\big)-\E\big(B^H_v-B^H_u\big)
\big(B^H_t-B^H_s\big)\big]\nonumber\\
&\quad=\frac 12\big[(t+u)^{2H}-(t+v)^{2H} +(s+v)^{2H}
-(s+u)^{2H}\big].
\end{align}
From the latter it follows that (note that, for $B^{H}$ the normalizing constant $\kap=1$, i.e., it coincides with that of sfBm)
\begin{equation}\label{e:kovar_skirtumas3}
d_{k,j}^{\widehat S^H,in}-d_{k,j}^{\widehat B^H,in}=\frac{\widehat\rho_H(j+k)}{4-2^{2H}}\quad\mbox{for}\ \vert j-k\vert\gs 1
\end{equation}
and that
\begin{equation}\label{e:kovar_skirtumas3a}
d^{\widehat S^H,in}_{kk}=1-\frac{b(k,H)}{4-2^{2H}}=d^{\widehat B^H,in}_{kk}-\frac{b(k,H)}{4-2^{2H}}\,,
\end{equation}
where
\[
b(k,H)=2^{2H-1}(k+1)^{2H}+3\cdot 2^{2H}k^{2H} +2^{2H-1}(k-1)^{2H}
-2(2k+1)^{2H}-2(2k-1)^{2H}.
\]
From  Lemma 1 \cite{km} (see also Lemma 6.24 \cite{KMR}) we get estimates
\begin{equation}\label{e:kovar}
\frac{\vert b(k,H)\vert}{4-2^{2H}}\ls \frac{33}{9\ln 4}\, k^{2H-4}\ls 3 k^{2H-4}\qquad\mbox{for}\ k\gs 3
\end{equation}
and
\begin{equation}\label{e:kovar1}
\max_{H\in(0,1)}  d^{\widehat S^H}_{kk}\ls \max_{H\in(0,1)}  d^{\widehat S^H}_{11}=\frac 76\,,\quad  \max_{H\in(0,1)}\frac{\vert b(k,H)\vert}{4-2^{2H}} \ls \max_{H\in(0,1)}\frac{\vert b(1,H)\vert}{4-2^{2H}}\ls \frac{1}{6}\,.
\end{equation}
Therefore,
\begin{align*}
(in)^{-1}\E V^{\widehat S^H}_{in,T}-1=& \frac{1}{in}\sum_{k=1}^{in-1}\big( d_{k,k}^{\widehat S^H,in} - d_{k,k}^{\widehat B^H,in}\big) -\frac{1}{in}\\
=& -\frac{1}{in}\sum_{k=1}^{in-1}\frac{b(k,H)}{4-2^{2H}}-\frac{1}{in} =O\left(\frac{1}{n}\right)
\end{align*}
and condition \eqref{e:dispersiju_skirtumas3} is satisfied. Further, by applying \eqref{e:kovar_skirtumas3}-\eqref{e:kovar}, and \eqref{e:vanishing0} we get
\begin{align*}
\max_{1\ls j\ls in-1}\sum_{k=1}^{in-1}\vert d_{k,j}^{\widehat S^H,in}\vert \ls& \max_{1\ls j\ls in-1}\sum_{k=1}^{in-1}\vert d_{k,j}^{\widehat S^H,in} -d_{k,j}^{\widehat B^H,in}\vert +\max_{1\ls j\ls in-1}\sum_{k=1}^{in-1}\vert d_{k,j}^{\widehat B^H,in}\vert\\
\ls&\max_{1\ls j\ls in-1}\sum_{k=1\atop k\neq j}^{in-1}\frac{\abs{\widehat\rho_H(j+k)}}{4-2^{2H}} +\max_{1\ls j\ls in-1}\vert d_{j,j}^{\widehat S^H,in}-d_{j,j}^{\widehat B^H,in}\vert\\
&+\max_{1\ls j\ls in-1}\sum_{k=1}^{in-1}\vert d_{k,j}^{\widehat B^H,in}\vert \ls C,
\end{align*}
where the last inequality is due to the fact that \eqref{e:maks_sum} applies to $B^H$ as well (see, e.g., Lemma 3.1 \cite{ks} and Lemma 2.7 \cite{KMR}). Hence \eqref{e:maks_sum} holds for sfBm. To finish the proof, we make use of Remark \ref{r:sufficientCond}. First we check the condition \eqref{e:kovariaciju_skirtumas1}. For the first term, by applying \eqref{e:kovar},  \eqref{e:kovar1} and $d^{\widehat B^H,n}_{k,k}=1$, we have
\begin{align*}
\frac 1n\sum_{k=1}^{n-1}\big\vert \big(d^{\widehat S^H,n}_{k,k}\big)^2 -\big(d^{\widehat B^H,n}_{k,k}\big)^2\big\vert=&
\frac 1n\sum_{k=1}^{n-1}\big\vert d^{\widehat S^H,n}_{k,k} -d^{\widehat B^H,n}_{k,k}\big\vert\cdot\vert d^{\widehat S^H,n}_{k,k}+d^{\widehat B^H,n}_{k,k}\big\vert\\
\ls&\frac C n\sum_{k=1}^{\infty}\frac{\abs{b(k,H)}}{4-2^{2H}}=O\left(\frac{1}{n}\right).
\end{align*}
For the second one, applying \eqref{e:kovar_skirtumas3}, we have
\begin{align*}
&\frac 1n\sum_{k=1}^{n-2}\sum_{m=1}^{n-k-1}\big\vert \big(d^{\widehat S^H,in}_{k,k+m}\big)^2 -\big(d^{\widehat B^H,in}_{k,k+m}\big)^2\big\vert\ls
\frac C n\sum_{k=1}^{n-2}\sum_{j=k+1}^{n-1}\frac{\abs{\widehat\rho_H(j+k)}}{4-2^{2H}}\\
&\quad\ls
\frac C n\sum_{k=1}^{n-2}\sum_{j=k+1}^{n-1}\frac{1}{j^{4-2H}}
\ls \frac C n\sum_{k=1}^{n-2}\int_k^{\infty}\frac{d x}{x^{4-2H}}\ls
\frac C n\sum_{k=1}^{\infty}{k^{2H-3}}=O\left(\frac{1}{n}\right)
\end{align*}
since $H\in(0,1)$.

Now we verify condition \eqref{e:kovariaciju_skirtumas2}. As above, we get
\begin{align}\label{e:isskaidymas1}
&\frac{1}{n} \sum_{j=1}^{n-1}\sum_{k=1}^{2n-1} \big\vert \big(c_{j,k}^{\widehat S^H,n}\big)^2 - \big(c_{j,k}^{\widehat B^H,n}\big)^2\big\vert
\ls \frac{C}{n} \sum_{j=1}^{n-1}\sum_{k=1}^{2n-1}\big\vert c_{j,k}^{\widehat S^H,n} - c_{j,k}^{\widehat B^H,n}\big\vert\,.
\end{align}
Let us denote  $\tau_n=\lfloor \sqrt n\rfloor$. Decompose  \eqref{e:isskaidymas1} in two sums in such a way
\begin{align*}
&\frac{1}{n} \sum_{j=1}^{n-1}\sum_{k=1}^{2n-1}\big\vert c_{j,k}^{\widehat S^H,n} - c_{j,k}^{\widehat B^H,n}\big\vert\\
&\quad\ls\frac{1}{n} \sum_{j=1}^{\tau_n}\sum_{k=1}^{2n-1}\big\vert c_{j,k}^{\widehat S^H,n} - c_{j,k}^{\widehat B^H,n}\big\vert +\frac{1}{n} \sum_{j=\tau_n+1}^{n-1}\sum_{k=1}^{2n-1}\big\vert c_{j,k}^{\widehat S^H,n} - c_{j,k}^{\widehat B^H,n}\big\vert=:J^{(1)}_n+J^{(2)}_n.
\end{align*}

From \eqref{e:difference} it follows that
\[
c_{j,k}^{\widehat S^H,n} - c_{j,k}^{\widehat B^H,n}=\frac{\widetilde\rho_H(2j+k)}{2^H(4-2^{2H})}\, .
\]
Applying the above equality and \eqref{e:vanishing} we obtain
\begin{align}\label{e:estimate}
J^{(1)}_n=&\frac{1}{n}\bigg(\frac{\vert\widetilde\rho_H(3)\vert}{2^H(4-2^{2H})}+ \sum_{j=2}^{\lfloor \sqrt n\rfloor} \frac{\vert\widetilde\rho_H(2j+1)\vert}{2^H(4-2^{2H})} \bigg) +\frac{1}{n} \sum_{j=1}^{\tau_n}\sum_{k=2}^{2n-1}\frac{\widetilde\rho_H(2j+k)}{2^H(4-2^{2H})}\nonumber\\
\ls& \frac{1}{n}\bigg(1+ C \sum_{j=2}^{\lfloor \sqrt n\rfloor}\frac{1}{(2j+1)^2}\bigg)+ \frac{C}{n} \sum_{j=1}^{\tau_n}\sum_{k=2}^{2n-1} \frac{1}{(2j+k)^2}\nonumber\\
\ls&\frac{1}{n}+ \frac{C}{\sqrt n}+\frac{C}{\sqrt n} \sum_{k=1}^\infty\frac{1}{k^2}=O\left(\frac{1}{\sqrt n}\right)
\end{align}
and
\begin{align*}
J^{(2)}_n\ls& \frac{C}{n} \sum_{j=\tau_n+1}^{n-1}\sum_{k=1}^{2n-1}\frac{1}{(2j+k)^{4-2H}}\ls \frac{C}{n} \sum_{j=\tau_n+1}^{n-1}\frac{1}{\tau_n+1}\sum_{k=1}^{2n-1}\frac{1}{k^{3-2H}}\\
\ls&   \frac{C}{\sqrt n}\sum_{k=1}^\infty\frac{1}{k^{3-2H}}=O\left(\frac 1{\sqrt n}\right).
\end{align*}
In \eqref{e:estimate} we used inequality
\[
\frac{\vert\widetilde\rho_H(3)\vert}{2^H(4-2^{2H})}\ls 1,
\]
which we obtained by computer simulation.

\section{Appendix}

Let us recall the condition $3 (e)$ of Theorem 2 in an article by B\'egyn \cite{begyn2} for the Gaussian process $(X_t)_{t\in[0,1]}$:

$3 (e)$ \emph{there exists a bounded function} $\widetilde g \dvit(0,1)\to R$ \emph{such that}
\[
\lim_{h\to 0+}\sup_{h\ls t\ls 1-2h}\bigg\vert \frac{(\delta^h_1\circ\delta^h_2 R)(t+h,t)}{h^{2\g}L(h)} -\widetilde g(t)\bigg\vert= 0,
\]
\emph{where} $R(s,t)$ \emph{is a covariance function of} $X$, $L:(0, 1)\to\mathbb{R}$ \emph{is a positive slowly varying function, $\g\in(0,1)$.}

For the Gaussian process $X$ we have
\begin{align*}
&(\delta^h_1\circ\delta^h_2 R)(t+h,t) \\ &\quad=\E\big(X_{t+2h}-2X_{t+h}+X_t\big)\big(X_{t+h}-2X_t+X_{t-h}\big)\\
&\quad= R(t+2h,t+h)-2R(t,t+2h)+R(t+2h,t-h)-2R(t+h,t+h)\\
&\qquad+5R(t,t+h)-2R(t+h,t-h)-2R(t,t)+R(t,t-h).
\end{align*}

1. Let $S^H=(S_t^H)_{t\in[0,1]}$ be a sfBm.  We prove that the condition $3 (e)$ of Theorem 2 in B\'egyn \cite{begyn2} does not hold for sfBm $S^H$ with $H\neq 1/2$.
Note that the following equality
\begin{align*}
\mu_t(h):=&\E\big(S^H_{t+2h}-2S^H_{t+h}+S^H_t\big)\big(S^H_{t+h}-2S^H_t+S^H_{t-h}\big)-\frac 12\big(2^{2H+2}-3^{2H}-7\big)h^{2H}\\
=&2^{2H+1}t^{2H}-\frac 12\,(2t-h)^{2H}-3(2t+h)^{2H}+2^{2H+1}(t+h)^{2H}-\frac 12\,(2t+3h)^{2H}
\end{align*}
holds. From the Taylor expansion for each fixed $t\neq mh$, $1\ls m\ls \lfloor h^{-1}-2\rfloor$, and small $h$ such that $(3h)/(2t)< 1$   we get
\begin{align*}
\mu_t(h)=&(2t)^{2H}\Big[2-\frac 12-3+2-\frac 12\Big]\\
&+(2t)^{2H}\sum_{k=1}^\infty \frac{2H(2H-1)\cdots(2H-k+1)}{k!}\Big(\frac{h}{2t}\Big)^k\Big[-\frac 12\,(-1)^k-3+2^{k+1}-\frac{3^k}{2}\Big] \\
=&(2t)^{2H}\sum_{k=4}^\infty \frac{2H(2H-1)\cdots(2H-k+1)}{k!}\Big(\frac{h}{2t}\Big)^k\Big[-\frac 12\,(-1)^k-3+2^{k+1}-\frac{3^k}{2}\Big].
\end{align*}
Since
\begin{align*}
&(2t)^{2H}\sum_{k=4}^\infty \frac{\abs{2H(2H-1)\cdots(2H-k+1)}}{k!}\,\bigg(\frac{ih}{2t}\bigg)^k \\
&\quad\ls (2t)^{2H}\sum_{k=4}^\infty \frac 2k\,\bigg(\frac{ih}{2t}\bigg)^k
\ls  - 2(2t)^{2H}\,\bigg(\frac{ih}{2t}\bigg)^3\ln \bigg(1-\frac{ih}{2t}\bigg)
\end{align*}
for $i=1,2,3$, then $\mu_t(h)=O(h^4)$. Thus,
\[
\frac{(\delta^h_1\circ\delta^h_2 G_H)(t+h,t)}{h^{2H}} \tends{h\to 0+} \frac 12 \big(2^{2H+2}-3^{2H}-7\big)\quad\mbox{for}\quad 0<H<1
\]
for $t\neq mh$, $1\ls m\ls \lfloor h^{-1}-2\rfloor$.
If $t= mh$, $1\ls m\ls \lfloor h^{-1}-2\rfloor$, $h\ls 1/3$, then using computer simulation we get
\begin{align*}
\frac{\mu_{mh}(h)}{h^{2H}}=&\frac{2^{2H+1}(mh)^{2H}- 2^{-1}(mh)^{2H}-3(3mh)^{2H}+2^{2H+1}(2mh)^{2H}-2^{-1}(5mh)^{2H}}{h^{2H}}\\
=& m^{2H}\big(2^{2H+1}-2^{-1}-3^{2H+1}+2^{2H+2}-2^{-1}5^{2H}\big)\neq 0.
\end{align*}
Thus, we haven't uniform convergence in condition $3 (e)$ for $H\in(0,1)$ and we can't apply the result of  B\'egyn \cite{begyn2} for the sfBm.

2. Let $B^{KH}=(B^{KH}_t)_{t\in[0,1]}$ be an fBm with index $KH$, and let $B^{H,K}=(B^{H,K}_t)_{t\in[0,1]}$ be a bifBm. Recall that the bifractional Brownian motions with parameters $H\in(0,1)$ and $K\in(0,1]$ is a centered Gaussian process with covariance function
\[
R_{HK}(t,s)=2^{-K}\big((t^{2H} +s^{2H})^K-\vert t-s\vert^{2HK}\big),\qquad s,t\gs 0.
\]
We prove that the condition $3 (e)$ of Theorem 2 in B\'egyn \cite{begyn2} does not hold for bifBm $B^{H,K}$ with $H\neq 1/2$ and $K\neq 1$.

Indeed, we can check that
\begin{align}\label{e:kovariacija}
\mu_t(h):=&\E\big(B^{H,K}_{t+2h}-2B^{H,K}_{t+h}+B^{H,K}_t\big)\big(B^{H,K}_{t+h} -2B^{H,K}_t+B^{H,K}_{t-h}\big) -2^{-K}\big(2^{2+2KH} -7-3^{2KH}\big)h^{2KH}\nonumber\\
=&2^{1-K}\big[2c(t,t+2h)-c(t+h,t+2h)-5c(t,t+h)-c(t-h,t+2h)\nonumber\\
&+2c(t-h,t+h)-c(t-h,t)\big] h^{2KH},
\end{align}
where
\[
c\big(t-ih,t+jh\big)=\frac{\s_{B^{H,K}}^2(t-ih,t+jh)}{2^{1-K} ((j+i)h)^{2KH}}- 1, \qquad i\in\{-1,0,1\},\ j\in\{0,1,2\},\quad i+j\neq 0.
\]
Since
\[
\s_{B^{H,K}}^2(t-ih,t+jh)=2^{1-K}\big( (i+j)^{2KH}h^{2KH}-f_{i,j,t}(h)\big)
\]
with
\[
f_{i,j,t}(h):=\big((t+jh)^{2H}+(t-ih)^{2H}\big)^K-2^{K-1}\big[(t+jh)^{2KH}+(t-ih)^{2KH}\big],
\]
we get that $f_{i,j,t}(0)=f^\prime_{i,j,t}(0)=0$. By Taylor's formula we obtain
\[
c\big(t-ih,t+jh\big)=-\frac{1}{{((i+j)h)^{2KH}}}\int_0^h f^{\prime\prime}_{i,j,t}(x)(h-x)\,dx,
\]
where
\begin{align}\label{e:bifBm_expansion}
f^{\prime\prime}_{i,j,t}(x)=& 4K(K-1)H^2\frac{\big[(t+jx)^{2H-1}j -(t-ix)^{2H-1}i \big]^2}{\big[(t+jx)^{2H}+(t-ix)^{2H}\big]^{2-K}}\nonumber\\
&+2HK(2H-1)\frac{(t+jx)^{2H-2}j^2 +(t-ix)^{2H-2}i^2}{\big[(t+jx)^{2H}+(t-ix)^{2H}\big]^{1-K}}\nonumber\\
&-2^K HK(2HK-1)\big[(t+jx)^{2KH-2}j^2 +(t-ix)^{2KH-2}i^2 \big].
\end{align}
Note that for $H> 1/2$, $j\in\{1,2\}$, and $t> x\gs 0$,
\begin{align*}
&\frac{[(t+jx)^{2H-1}j -(t-ix)^{2H-1}i \big]^2}{[(t+jx)^{2H}+(t-ix)^{2H}]^{2-K}} \ls\frac{(i^2+j^2)(t+jx)^{4H-2}}{[(t+jx)^{2H} +(t-ix)^{2H}]^{2-K}} \\
&\quad=(i^2+j^2)\bigg[\frac{(t+jx)^{2H}}{(t+jx)^{2H} +(t-ix)^{2H}}
\bigg]^{2-K}\,(t+jx)^{2HK-2}\\
&\quad\ls \frac{i^2+j^2}{(t+jx)^{2-2KH}}\ls \frac{i^2+j^2}{(t-x)^{2-2KH}}\,.
\end{align*}
If $j=0$, then for the first term of \eqref{e:bifBm_expansion}, we get
\[
\frac{i^2(t-ix)^{4H-2}}{[t^{2H}+(t-ix)^{2H}]^{2-K}}\ls \frac{1}{(t-x)^{2-2KH}}\,.
\]
If $H< 1/2$ and $t>x\gs 0$, then
\[
\frac{[(t+jx)^{2H-1}j -(t-ix)^{2H-1}i \big]^2}{[(t+jx)^{2H}+(t-ix)^{2H}]^{2-K}} \ls\frac{(i^2+j^2)(t-ix)^{4H-2}}{(t-ix)^{4H-2KH}} =\frac{i^2+j^2}{(t-x)^{2-2KH}}\,.
\]
Thus, for $t>h\gs 0$,
\begin{align*}
\sup_{0\ls x\ls h}\vert f^{\prime\prime}_{i,j,t}(x)\vert\ls& \frac{4(i^2+j^2)}{(t-h)^{2-2KH}}\,{\bf 1}_{\{H\gs 1/2\}} +\frac{4(i^2+j^2)}{(t-h)^{2-2KH}}\,{\bf 1}_{\{H< 1/2\}}\\
&+\frac{2(i^2+j^2)}{[(t-h)^{2H}]^{1-K}(t-h)^{2-2H}}+\frac{2(i^2+j^2)}{(t-h)^{2-2KH}}
= \frac{8(i^2+j^2)}{(t-h)^{2-2KH}}
\end{align*}
and
\begin{align*}
&\sup_{0<\d\ls h}\big\vert c\big(t-ih,t+jh\big)\big\vert\ls \sup_{0<\d\ls h}\frac{3(i^2+j^2) h^{2-2KH}}{(i+j)^{2KH}(t-h)^{2-2KH}}\tends{h\to 0+}0.
\end{align*}
Therefore, for each fix $t>h$, it follows that
\[
\frac{\mu_t(h)}{h^{2KH}}\tends{h\to 0+}0.
\]
Now we return to equality \eqref{e:kovariacija}. After simple calculation for $0\ls u<v\ls s<t$, we get
\begin{align*}
&\E\big[\big(B^{H,K}_v-B^{H,K}_u\big)\big(B^{H,K}_t-B^{H,K}_s\big) -2^{1-K}\E\big(B^{KH}_v-B^{KH}\big)\big(B^{KH}_t-B^{KH}_s\big)\big]\\
&\quad=\frac{1}{2^K}\big[(t^{2H}+v^{2H})^K-(t^{2H}+u^{2H})^K -(s^{2H}+v^{2H})^K+(s^{2H}+u^{2H})^K\big].
\end{align*}
Thus
\begin{align*}
\mu_t(h)
=&2^{-K}\big[ \big((t+2h)^{2H}+(t+h)^{2H}\big)^K
-2\big((t+2h)^{2H}+t^{2H}\big)^K-(2^K+1)(t+h)^{2KH}\\
&+(2^{1-K}+3)\big((t+h)^{2H}+t^{2H}\big)^K -(2^K+1)t^{2KH}+\big((t+2h)^{2H}+(t-h)^{2H}\big)^K \\
& -2\big((t+h)^{2H}+(t-h)^{2H}\big)^K+\big(t^{2H}+(t-h)^{2H}\big)^K\big].
\end{align*}
Then for $K\neq 1$ and a sequence of positive numbers $(h_n)$  converging to zero, using computer simulation, we obtain
\begin{align*}
\frac{\mu_{h_n}(h_n)}{h_n^{2KH}}=&2^{-K}\big[ \big(3^{2H}+2^{2H}\big)^K
-2\big(3^{2H}+1\big)^K-(2^K+1)2^{2KH} \\
&+(2^{1-K}+3)\big(2^{2H}+1\big)^K-(2^K+1) +3^{2KH} -2^{2KH+1}+1\big]\neq 0.
\end{align*}
So we have no uniform convergence in condition $3 (e)$ for $H,K\in(0,1)$, $H\neq 1/2$,  and thus we cannot apply the result of  B\'egyn \cite{begyn2} for the bifBm.


\begin{thebibliography}{00}

\bibitem{BS11}
J.-M. Bardet, D. Surgailis,  Measuring the roughness of random paths by
  increment ratios, \emph{Bernoulli}, \textbf{7}~(2) (2011), 749--780.

\bibitem{BS13}
J.-M. Bardet, D. Surgailis, Moment bounds and central limit theorems
  for {G}aussian subordinated arrays, \emph{Journal of Multivariate Analysis}, \textbf{114} (2013),
  457--473.

\bibitem{baxter} G. Baxter, A strong limit theorem for Gaussian processes, \emph{Proc. Am. Math. Soc.} \textbf{7}, (1956), 522--527.


\bibitem{begyn1} A. B\'egyn, Quadratic variations along irregular partitions for Gaussian processes, \emph{Electronic Journal of Probability}, \textbf{10} (2005), 691-717.

\bibitem{begyn2} A. B\'egyn, Asymptotic development and central limit theorem for
quadratic variations of Gaussian processes, \emph{Bernoulli} \textbf{13}(3) (2007), 712-753.


\bibitem{BGT} T. Bojdecki, L. G. Gorostiza, A. Talarczyk, Sub-fractional Brownian motion and its relation to occupation time, Stat. \& Probab. Lett. 69 (2004), 405-419.




\bibitem{Coeurjolly-01}  J.-F. Coeurjolly, Estimating the parameters of a fractional Brownian motion by discrete variations of its sample paths, \emph{Statistical Inference for Stochastic Processes}, \textbf{4} (2001), 199--227.

\bibitem{cohen} S. Cohen, X. Guyon, O. Perrin, and M. Pontier,  Singularity functions for fractional processes, and application to fractional Brownian sheet, \emph{Ann. Inst. Henri Poincar\'e, Probab. Stat.}, \textbf{42} (2005), 187-205.

\bibitem{glad}  E. G. Gladyshev, A new Limit theorem for stochastic processes with Gaussian increments, \emph{Theory Probab. Appl.}, \textbf{6}(1) (1961),  52-61.

\bibitem{guyon} X. Guyon and J. Leon, Convergence en loi des H-variations d'un processus gaussien stationnaire sur R, \emph{Ann. Inst. Henri Poincar\'e, Sect. B, Probab. Stat.}, \textbf{25}(3) (1989), 265-282.

\bibitem{HV} C. Houdr\'e and J. Villa, An example of infinite dimensional quasi-helix, Contemporary Mathematics 366 (2003), 195-201.

\bibitem{Iss} L. Isserlis, On a formula for the product-moment coefficient of any order of a normal frequency distribution in any number of variables. \emph{Biometrika} \textbf{12} (1918) 134-139.

\bibitem{IL} J. Istas and G. Lang, Quadratic variations and estimation of the local H\"older index of a Gaussian process. \emph{Ann. Inst. Henri Poincar\'e, Probab. Stat.}, \textbf{33} (1997),   407-436.

\bibitem{Istas98}
A. Benassi, S. Cohen, J. Istas and S. Jaffard, Identification of filtered white noises, \emph{Stochastic Processes and their Applications}, \textbf{75}~(1) (1998), 31--49.

\bibitem{kg} R. Klein and E. Gin\'e. On quadratic variation of processes with Gaussian
increments. \emph{Ann. Probab.}, \textbf{3}(4) (1975)  716-721.

\bibitem{Kolmogorov40}
A. Kolmogorov, Wienersche {S}piralen und einige andere interessante
  {K}urven im {H}ilbertschen raum. \emph{Proceedings of the USSR Academy of Sciences},
  \textbf{26} (1940) 115--118.

\bibitem{kk1} K. Kubilius, On estimation of the Orey index for a class of Gaussian processes, Stochastics: An International Journal Of Probability And Stochastic Processes, \textbf{87}(4) (2015) 562-591.

\bibitem{kubmel1} K. Kubilius, D. Melichov, Quadratic variations and estimation of the Hurst index of the solution of SDE driven by a fractional Brownian Motion, \emph{ Lithuanian Mathematical Journal}, \textbf{50}(4)  (2010)  401-417.

\bibitem{km}  K. Kubilius, D. Melichov, Exact confidence intervals of the extended Orey index for Gaussian processes, Methodol. Comput. Appl., \textbf{18}(3), (2016) 785--804.

\bibitem{ks} Kubilius, K. and Skorniakov, V., On some estimators of the {H}urst index of the solution of {SDE} driven by a fractional {B}rownian motion, Statist. Probab. Lett., \textbf{109}, (2016) 159--167.

\bibitem{ksm} K. Kubilius, V. Skorniakov, and D. Melichov, Estimation of parameters of SDE driven by fractional Brownian motion with polynomial drift, \emph{J. Stat. Comput. Simul.} \textbf{86}(10)  (2016) 1954-1969.

\bibitem{KMR} K. Kubilius, Yu. Mishura, K. Ralchenko, Parameter Estimation in Fractional Diffusion Models, Springer (Bocconi \& Springer Series), 2017.


\bibitem{ma} R. Malukas, Limit theorems for a quadratic variation of Gaussian processes. \emph{Nonlinear Analysis: Modelling and Control}, \textbf{16}(4) (2011), 435-452.

\bibitem{Mandelbrot68}
B.~B. Mandelbrot, J. W.~V. Ness,  Fractional {B}rownian motions,
  fractional noises and applications. \emph{SIAM Review} \textbf{10}(4) (1968) 422--437.

\bibitem{rn0} R. Norvai\v sa and D.M. Salopek,  Estimating the Orey Index of a
Gaussian Stochastic Process with Stationary Increments: An Application
to Financial Data Set. Canadian Mathematical Society Conference Proceedings,
\textbf{26}  (2000), 353-374.

\bibitem{rn1} R. Norvai\v sa, A coplement to Gladyshev's theorem. \emph{Lithuanian Mathematical Journal}, \textbf{51}(1)  (2011), 26-35.

\bibitem{rn2} R. Norvai\v sa, Weighted power variation of integrals with respect to a Gaussian process. \emph{Bernoulli}, \textbf{21}(2) (2015) 1260-1288.


\bibitem{or} S. Orey,  Gaussian sample functions and the Hausdorff dimension of level crossings. \emph{Z. Wahrsch. verw. Gebiete}, \textbf{15} (1970)  249-256.


\bibitem{viit} L. Viitasaari, Necessary and sufficient conditions for limit theorems for quadratic variations of Gaussian sequences. \emph{Probability Surveys}, \textbf{16} (2019) 62-98.


\end{thebibliography}
\end{document}